\documentclass[11pt,a4paper,reqno]{amsart}

\usepackage{amsfonts,amsthm,amscd,epsfig,amsmath,amssymb,enumerate}

\author {A. Faggionato}
\address{Dipartimento di Matematica ``G. Castelnuovo", Universit\`a   ``La
  Sapienza'', P.le Aldo Moro  2, 00185  Roma, Italy.}
\email{faggiona@mat.uniroma1.it}

\numberwithin{equation}{section}

\DeclareMathSymbol{\leqslant}{\mathalpha}{AMSa}{"36} % nicer `smaller or equal'
\DeclareMathSymbol{\geqslant}{\mathalpha}{AMSa}{"3E} % nicer `larger or equal'
\DeclareMathSymbol{\eset}{\mathalpha}{AMSb}{"3F}     % nicer `emptyset'
\renewcommand{\leq}{\;\leqslant\;}                   % redef. of < or =
\renewcommand{\geq}{\;\geqslant\;}                   % redef. of > or =
\usepackage{graphicx}
\usepackage{url}

\newtheorem{theorem}{Theorem}

\newtheorem{lemma}[theorem]{Lemma}

\newtheorem{prop}[theorem]{Proposition}

\newtheorem{theorem*}{Theorem}

\numberwithin{theorem}{section} 

\newcommand{\ignore}[1]{}
\newcommand{\be}{\begin{equation}}

\def\1{\ifmmode {1\hskip -3pt \rm{I}} \else {\hbox {$1\hskip -3pt \rm{I}$}}\fi}

\newcommand{\cA}{\ensuremath{\mathcal A}}
\newcommand{\cB}{\ensuremath{\mathcal B}}
\newcommand{\cC}{\ensuremath{\mathcal C}}

\newcommand{\cF}{\ensuremath{\mathcal F}}
\newcommand{\cG}{\ensuremath{\mathcal G}}

\newcommand{\cL}{\ensuremath{\mathcal L}}

\newcommand{\bbE}{{\ensuremath{\mathbb E}} }

\newcommand{\bbL}{{\ensuremath{\mathbb L}} }

\newcommand{\bbN}{{\ensuremath{\mathbb N}} }

\newcommand{\bbP}{{\ensuremath{\mathbb P}} }

\newcommand{\bbZ}{{\ensuremath{\mathbb Z}} }

    \let\d=\delta  \let\e=\varepsilon
  \let\h=\eta      
      \let\o=\omega      
  \let\s=\sigma    
  
\let\D=\Delta      
\let\O=\Omega

\newcommand{\ZZ}{{\mathbb Z}}

\newcommand{\Aa}{{\mathcal A}}

\newcommand{\Bb}{{\mathcal B}}

\newcommand{\Ff}{{\mathcal F}}

\newcommand{\II}{{\mathbb I}}
    \let\d=\delta  \let\e=\varepsilon
  \let\h=\eta

\begin{document}

\title[%Quenched CLT %invariance principle
[Hydrodynamic limit of symmetric exclusion processes]{Hydrodynamic
limit of symmetric  exclusion processes in inhomogeneous  media}

\begin{abstract}
In \cite{J} M. Jara has presented a method,   reducing the proof of
the hydrodynamic limit of symmetric exclusion processes  to an
homogenization problem, as unified approach to recent works on the
field as \cite{N}, \cite{F1}, \cite{F2} and \cite{FJL}. Although not
stated in \cite{J}, the reduction of the hydrodynamic limit to an
homogenization problem was already obtained  (in a different way) in
\cite{N}, \cite{F1}. This alternative and very simple relation
between  the two problems  goes back to an idea of K.\ Nagy
\cite{N},   is
 stated in \cite{F1}[Section B] for exclusion processes on
$\bbZ^d$ and, as stressed in \cite{F2}, is completely general. The
above relation has been applied to \cite{N}, \cite{F1}, \cite{F2}
and \cite{FJL}
  and could be applied to  other symmetric exclusion processes,
mentioned in \cite{J}. In this short note we briefly recall this
unified approach  in a complete general setting. Finally, we recall
how the homogenization problem has been solved in the above previous
works.

\medskip

\noindent \textsl{Key words}: exclusion process, homogenization,
disordered system, hydrodynamic limit.

%\medskip

%\noindent \textsl{AMS 2000 subject classification}:
%
%60K37, % Processes in random environments
%
%60F17, % Functional limit theorems; invariance principles
%
%60G55. %Point processes
\end{abstract}

\maketitle

\section{Introduction}

Exclusion processes in inhomogeneous media can be used to analyze
transport properties of particle systems in inhomogeneous and/or
disordered media with hard--core interaction. Due to the possible
lack of spacial invariance, the resulting model is an example of
non--gradient exclusion process, in the sense that the transition
rates cannot be written as gradient of some local observable.
 Despite this fact, if the jump rates do not depend on the orientation of the jump
(as in \cite{N}, \cite{F1}, \cite{F2}, \cite{FJL} but not as in
\cite{Q1}, \cite{Q2} and \cite{FM}),
  the hydrodynamic
limit of the exclusion process can be proven without using the very
sophisticated techniques developed for non--gradient systems  (cf.
\cite{KL} and references therein), which in addition would require
non trivial spectral gap estimates  that fail in the case of jump
rates  non bounded from below by a positive constant (cf. Section
1.5 in \cite{M}). The strong simplification comes from the fact that
the infinitesimal variation of the occupancy  number at a give site
is a linear combination of occupancy numbers.

Alternative routes to prove the hydrodynamic limit are the method of
{\sl corrected empirical measures} developed in \cite{JL} and
\cite{GJ}, which is also at the base of \cite{J}, and a duality
relation developed in \cite{N}, \cite{F1}. The advantage of the
former is that it  works also for  zero range processes, while the
latter is very elementary and direct. In both cases, the
hydrodynamic limit reduces  to an homogenization result for a random
walk (possibly in a random environment). In \cite{J} the connection
between hydrodynamics and homogenization by means of the corrected
empirical measure has been described in a general setting, and
presented as ``unified approach" to recent works on the field as
\cite{N}, \cite{F1}, \cite{F2} and \cite{FJL}. The alternative
connection due to the duality relation presented in \cite{N} and
\cite{F1} is not mentioned at all in \cite{J}. On the other hand,
this is indeed the connection  between hydrodynamics and
homogenization used in \cite{N}, \cite{F1}, \cite{F2} and
\cite{FJL}, and could be applied as well to the other symmetric
exclusion processes. In what follows we briefly recall this
equivalence working in a more general setting than \cite{N},
\cite{F1}.

Finally, we remark that, having  established the equivalence between
the two problems, it remains to solve the homogenization problem. In
\cite{N}  the author proves indeed an invariance principle for the
random walk, under rather restrictive conditions. In \cite{F1} the
author applies the method of Stone \cite{S} for random walks that
can be written  as space--time change of  Brownian motion, solving
the hydrodynamic limit for symmetric exclusion processes on $\bbZ$
under very natural conditions. Again, Stone's theory has been
applied in \cite{FJL}. In \cite{F2} the homogenization problem has
been solved by means of the notion of {\em two--scales convergence}
(\cite{Nu}, \cite{A}, \cite{ZP}) and is indeed the core of
\cite{F2}. In particular,  in \cite{J} the author does not  recover
the previous results of \cite{N}, \cite{F1}, \cite{F2} and
\cite{FJL} as stated in the Introduction of \cite{J}, since  the
homogenization problem is not solved in \cite{J}.

\section{Models and main result}

For each $n \geq 1$ let $G_n=(V_n, E_n)$ be a graph with unoriented
bonds. $V_n$ and $E_n$  denote respectively  the family of vertexes
and the family of edges.
 We suppose that $V_n $ is a locally finite subset of $X$, where   $X$ is  a metric space endowed of a
Radon measure $m$ (one could as well simply require that $V_n$ can
be injected in $X$ with locally finite image, this would slightly
change the notation below). Moreover, we assume that
\begin{equation}\label{vagamente}m_n:= \frac{1}{a_n} \sum_{x \in V_n} \d_x \rightarrow m
\end{equation} with respect to the vague topology, for a fixed sequence
$a_n$ such that
\begin{equation}\label{diva}
\lim_{n \to \infty} a_n = \infty\,.
\end{equation}

The state space of the exclusion process on $G_n$ is given by the
set $\{0,1\} ^{V_n}$, endowed of the product topology.
%Given a measurable nonnegative function $\rho$ on $X$, we say
%that the sequence $(\mu_n)_{n \geq 1}$, $\mu_n$ being a probability
%measure on $\{0,1\}^{V_n}$, is associated to the macroscopic profile
%$\rho$ if, under $\mu_n$, the random measures $\p_n$ converge
%vaguely to the measure $\rho(u)m(du)$ on $X$. This means that, given
%$G \in C_c (X)$ and $\d>0$, it holds
%\begin{equation}
%\lim _{n\uparrow \infty} \mu_n \Big( \Big| \frac{1}{a_n} \sum_{x \in
%V_n} G(x) \eta(x) - \int G(u) \rho(u) m (du)\Big|> \d \Big)=0\,.
%\end{equation}

The exclusion process on the graph $G_n$ is defined in terms of a
rescaling constant $b_n$ and a conductance field $\{ c_n(b) : b \in
E_n \}$. It is a Markov process on $\{0,1\}^{V_n}$ whose
  Markov generator acts  on local
functions (i.e. depending only on $\eta_x$ with $x$ varying in a
fixed finite subset of $V_n$) as
\begin{equation}
\cL_n f (\eta)= b_n \sum _{\{x,y\}\in E_n } c_n (x,y) \bigl( f
(\eta^{x,y})- f(\eta) \bigr) \,.
\end{equation}
Above we have written $c(x,y)$ for the conductivity $c ( \{x,y\})$
associated to the bond $\{x,y\}$ and we have used the standard
notation
$$
\eta^{x,y} _z=
 \begin{cases} \eta_y & \text{ if } z=x\,,\\
 \eta_x & \text{ if  } z=y \,,\\
 \eta_z & \text{ otherwise}\,.
 \end{cases}
 $$
We assume that the above dynamics is well defined (cf. \cite{L} for
sufficient conditions). It can be constructed as follows. On a
common probability space $(\O, \cF, \bbP)$ define a family of
independent Poisson processes $N^n_b(\cdot)$, parameterized by $b
\in E_n$, such that $\bbE \bigl( N^n_b (t) \bigr)= c_n (b) t$.
Roughly, if $t$ is a jump time of $N_b^n$ then at time $t$ the
exclusion process performs the jump $ \eta(t-)\mapsto \eta (t)^{b}$.
To formalize the above definition we assume  that  there exists
$\e>0$ such that  for almost all $\o$   the connected components of
the graph $\cG^n_\e (\o)$, obtained from $G_n$ by keeping only the
bonds $b$ such that $N^n_b (\e) [\o]\geq 1$, have finite cardinality
(cf. \cite{D}).
% % Form now on we assume this two properties.

In what follows, given $t \geq 0$ we denote by $\cG^n _t (\o)$ the
random graph with vertexes $V_n$ and edges  $b\in E_n $ such that
$N^n_b (t)[\o]> N^n _n (k\e )[\o]$, where $k\e \leq t < (k+1)\e$, $k
\in \bbN$. As consequence of the above assumption, there exists a
measurable subset $\cA\in \cF $, with $\bbP(\cA)=1$,  such that
$\cG^n_t(\o)$ has only connected components with finite cardinality
for all $t\geq 0$ and all $\o \in \cA$. Let $\o\in \cA $. Then,
given an initial configuration $\h(0)$, the configuration $\h (t)=
\h (t)[\o]$ at time $t$ is defined  as follows:

   Let $\cC$ be any connected
component of $ \cG^n_t(\o)$ and let
$$
\{s_1<s_2<\cdots <s_r\}= \{s\,:\, N^n_b(s)= N^n_b(s-)+1 ,\; b \text{
bond in } \cC ,\; 0<s\leq t\}.
 $$
 Start with $\h(0)$. At time $s_1$ switch the values between $\h_x$
 and $\h_{y}$ if $b=\{x,y\}$, $N^n_b (s_1)=N^n_b(s_1-)+1$ and $b$ is a bond in $\cC$.
 Repeat the same operation orderly for times $s_2,s_3, \dots , s_r$.
Then the  resulting configuration coincides with $\h(t)$ on $\cC$.

Below,  we denote by $\bbP_{\mu_n}$ the  law of the exclusion
process on $G_n$ when the initial distribution is $\mu_n$.
% We point
%out that if for example the graph $G_n$  has finite degree and $\sup
%_{b \in E_n }c_n(b) < \infty$, then the dynamics is always well defined.
%For simplicity, below we assume that  $\sup _{b \in E_n}c_n(b) < \infty$.

\medskip

As special case one can consider the dynamics of a single particle.
The resulting process is a continuous--time random walk on $G_n$,
having $V_n$ as state space and  infinitesimal generator
\begin{equation} \bbL_n h(x)= b_n \sum _{y\in V_n: y \sim x} c_n (x,y)
\bigl(h(y)- h(x) \bigr)\,.
\end{equation}
%defined on the family of bounded functions on $V_n$. Above we have
%used the convention:
Note that, given vertexes $a,b \in V_n$, we write $a \sim b $ if
$\{a,b\} \in E_n$. We write $X_n(t|x)$ for the above random walk
when starting at site $x$.

%A fundamental identity relating the exclusion process and the random
%walk if the following
%$$ \bigl(\cL_n \pi_n ( \eta)\bigr) (x)= (\bbL_n \eta) (x)$$

Given a function $\varphi$ on $V_n$, we write
 $$P_t^n \varphi (x):= \bbE\bigl[ \varphi (\,X_n(t|x)\,) \bigr]\,,
 \qquad x \in V_n\,.$$
%where $E [\cdot]$ denotes the expectation.

\begin{theorem}\label{annibale}
Suppose that there exists a stochastic process $W$ on $X$, with
symmetric probability kernel,  satisfying the following
homogenization property:
\begin{equation}\label{hom1}
\lim _{n\uparrow \infty} \frac{1}{a_n} \sum_{x \in V_n} \bigl| P_t^n
\varphi (x)- P_t \varphi(x) \bigr| =0 \,, \qquad \forall \varphi \in
C_c (X) \,,
\end{equation}
where
$$ P_t \varphi (x) := \bbE\bigl[ \varphi (\,W(t|x)\,) \bigr]$$ and $W(t|x)$ denotes the process $W$ starting at $x\in
X$. Then the following hydrodynamic behavior of the exclusion
process holds:
\smallskip

Let $\rho_0: X\to [0, \infty)$ be a bounded Borel function and let
$(\mu_n)_{n \geq 1}$ be a sequence of probability measures on
$\{0,1\}^{V_n}$ such that, for any $\d>0$ and any $\varphi \in \{
P_t\Psi: t \geq 0, \psi \in  C_c(X)\} $, it holds
\begin{equation}\label{prof1}
\lim _{n\uparrow \infty} \mu_n \Big( \Big| \frac{1}{a_n} \sum_{x \in
V_n} \varphi(x) \eta(x) - \int \varphi (u) \rho_0(u) m (du)\Big|> \d
\Big)=0\,.
\end{equation}
Then, for any $\d>0$, $t>0$ and  and any $\varphi  \in C_c (X)$ it
holds
\begin{equation}\label{prof2}
\lim _{n \uparrow \infty} \bbP_{\mu_n} \Big( \Big| \frac{1}{a_n}
\sum _{x \in V_n} \varphi (x) \eta ^n _t (x)- \int \varphi(u) (P_t
\rho_0) (u) m (du) \Big|> \d \Big)=0\,.
\end{equation}
\end{theorem}
We point out that, in applications, the set $\{ P_t\Psi: t \geq 0,
\psi \in C_c(X)\} $ is given by enough regular functions decaying to
$0$ at infinity.

One can restate the homogenization property \eqref{hom1} in a form
more related to partial differential equations as done in \cite{F2},
by means of the concept of 2--scale convergence which is particular
suited for the above setting (cf. \cite{ZP}). \eqref{hom1} is indeed
the starting point of \cite{F2}.

\smallskip

We point out  that the above structures, given by the   graphs
$G_n$,  the conductance fields $c_n$ and the initial distributions
$\mu_n$, could be random.  In order to obtain for example a quenched
hydrodynamic limit, it is enough to prove \eqref{hom1} for almost
all realizations of the {\em random environment}, given by the graph
$G_n$ and  the conductance field $c_n$.

\section{Proof of Theorem \ref{annibale}}

We write $p_n(t,x,y)$ for the probability that the random walk $X_n
(t|x)$ is in $y$ at time $t$. By the symmetry of the jump rates it
holds
\begin{equation}\label{simone}
p_n (t,x,y)= p_n (t,y,x)\, , \qquad \forall x, y \in V_n\,.
\end{equation}

Recall the graphical construction of the simple exclusion process
given in the previous section and set $N^n_{x,y}$ for
$N^n_{\{x,y\}}$. Since
$$d \eta_x (t) = \sum_{y \in V_n : y \sim x } ( \eta_y- \eta_x) (t-)
d N^n_{x,y} (t)$$ we can write
\begin{equation}\label{derivo}d \eta (t) = \bbL_n \eta (t) dt + dM^n(t), \end{equation} where
$$
dM^n_x (t)= \sum_{y \in V_n : y \sim x } ( \eta_y- \eta_x) (t-) d
A^n_{x,y} (t)\,,\qquad A^n_{x,y}(t)= N_{x,y}^n (t)- c(x,y) t \,.
$$
Note that
 $M^n_x(\cdot)$ has  trajectories of bounded
variation on finite intervals a.s.

 Formally, (\ref{derivo}) implies that
\begin{equation}\label{formale}
\h (t) = T(t) \h (0) + \int_0 ^t T(t-s) dM^n(s)
\end{equation}
where $T(t) = e^{t\bbL_n}$, i.e.
\begin{equation}\label{formale1}
\h _x (t) = \sum_{y \in V_n} p_n(t,x,y) \h _y (0) + \sum _{y \in
V_n} \int _0 ^t p(t-s,x,y) dM^n_y(s).
\end{equation}
Due to the graphical construction of the dynamics, if $\sum _{x\in
V_n}\h_x(0)<\infty$, then for all but a finite family of indexes $y$
 $ dM^n_y(s)=0$ for all $0\leq s\leq t$ and in particular  the last series
in (\ref{formale1}) reduces to a finite sum and is  meaningful. In
this case, one can check that
 (\ref{formale1}) holds a.s. by direct computation
 using that
$$
\frac{d}{dt} p_n(t, x,y)= \left( \bbL_n  p_n(t, \cdot, y)\right)_x.
$$
The following result shows that the site exclusion constraint is
negligible from a  hydrodynamic viewpoint:
\begin{prop}\label{gaetano} Given  $\d>0$, $t>0$, $\varphi \in C_c
(X)$ and given a sequence of probability measures $\mu_n$ on
$\{0,1\}^{V_n}$, it holds
\begin{equation}\label{limitone}
\lim _{n\uparrow\infty} \bbP _{\mu_n} \Big( \Big| \frac{1}{a_n} \sum
_{x\in V_n } \varphi(x) \h _x (t)- \frac{1}{a_n} \sum _{x\in V_n}
\varphi(x)\sum _{y \in V_n } p_n(t,x,y) \h _y (0) \Big|>\d \Big)=0.
\end{equation}
\end{prop}
\begin{proof}
 Let the support of $\varphi$ be included in $\D$ and fix $\e>0$.
Given $x\in V_n $ and $t>0$ define $\cC ^n  _x (t)$ as the connected
component of $\cG^n  _t$ containing $x$.  Then for each positive
integer $n $ we can choose a bounded set $B_n \subset X$ such that $
P(\Aa _n ^c ) <\epsilon$ where $\Aa_n$ is the subset of
configurations $\o$ satisfying the
following conditions:\\
\begin{align}
& \cup _{x\in V_n \cap \D } \cC^n_x (t) \subset B_n,\label{roma1}\\
& \frac{  \|\varphi \|_\infty}{a_n} \sum _{x \in V_n \cap \D  } \sum
_{y\in V_n \setminus B_n   } p_n (t,x,y ) \leq \d/2 .\label{roma2}
\end{align}
Given $\h (0)$ and $n$, we define $\h^{(n)} (0) \in \{0,1\}^\ZZ$ as
$ \h^{(n)} _x ( 0) = \h _x (0)\II_{x \in V_n \cap B_n} $ and write
$\h ^{(n)}(s)$ for  the configuration at time $s$ obtained by the
graphical construction when starting from $\h^{(n)}(0)$ at time $0$.

 Due to the graphical construction of the dynamics and
condition (\ref{roma1}), if $\o\in \Aa_n$ then $ \h _x  (t ) = \h _x
^{(n)} (t)$ for all $x \in V_n \cap B_n$, thus implying
$$
 \frac{1}{a_n} \sum
_{x \in V_n } \varphi(x)\h _x (t)=
 \frac{1}{a_n} \sum
_{x \in V_n } \varphi(x) \h _x ^{(n)} (t).
$$
 Moreover, due to
(\ref{roma2}), if $\o\in \Aa_n$ then
$$
 \frac{1}{a_n} \sum
_{x\in V_n} \bigl|\varphi(x)\bigr|\sum _{y \in V_n } p_n(t,x,y) |\h
_y(0)- \h _y ^{(n)} (0)|\leq \d/2. $$ Therefore the l.h.s. of
(\ref{limitone}) with fixed $n$ can be bounded by
$$ \bbP_{\mu_n} \left(  | Z_n |>\d /2 \right)+ P ( \Aa ^c_n ) \leq
4 \bbE_{\mu_n} \left( Z_n ^2 \right)/\d^2 +\e\,,
$$
 where
$$
Z_n =
 \frac{1}{a_n} \sum
_{x \in V_n } \varphi(x) \h ^{(n)} _x (t)- \frac{1}{a_n} \sum _{x\in
V_n } \varphi(x)\sum _{y \in V_n } p_n(t,x,y) \h _y^{(n)} (0).
$$
Since $\sum _{x\in V_n} \h ^{(n)} _x (0) <\infty$, setting
$$
dM^{(n)} _x (s)=\sum _{y \in V_n: y \sim x} (\h^{(n)} _{y}-\h^{(x)}
_k)(s-) d A^n_{x,y} (s),
$$
 $(\ref{formale1}) $ implies that
 \begin{equation*}
% \begin{split}
Z_n  =  \frac{1}{a_n} \sum _{x\in V_n } \varphi(x) \sum _{y \in V_n}
\int _0 ^{t} p_n(t -s,x,y) dM^{(n)}_y (s).
% & = \frac{1}{n} \sum _{k\in \ZZ} \varphi\left(\frac{k}{n}\right) \sum
%_{j\in \ZZ} \int _0 ^{tn^2}\left( p (t n^2-s,k, j) - p(t-s,k,
%j+1)\right) (\h ^{(n)} _{j+1} - \h ^{(n)} _j ) (s- ) d A _j (s)
%\end{split}
\end{equation*}
In order to conclude the proof it is enough to apply Lemma
\ref{mezzo} below to the above estimates.
\end{proof}

\begin{lemma}\label{mezzo}
For each $n\geq 1$ let $\nu _n $ be a probability measure on
$\{0,1\}^{V_n}$ such that $$ \nu _n \bigl(\sum _{x \in V_n } \h_x
<\infty\bigr) = 1\,. $$ Then
\begin{equation*}
 \lim _{n\uparrow\infty}
\bbE _{\nu _n} \Big[\Big(
\frac{1}{a_n} \sum _{x \in V_n } \varphi(x) \sum _{y\in V_n} \int _0
^{t} p_n(t-s,x,y) dM^n _y (s)
\Big)^2 \Big]=0.
\end{equation*}
\end{lemma}
Recall that the above series over $y$ reduces to a finite sum
whenever $\sum _{x \in V_n} \h_x(0) <\infty$, and therefore it is
well defined a.s.
\begin{proof}
We define $f_n$ as
\begin{multline}\label{mozart}
f_n= \frac{1}{a_n} \sum _{x\in V_n} \varphi(x)
 \sum _{y \in V_n}  \int _0 ^{t } p_n(t -s,x,y) dM ^n_y (s) =\\
 \frac{1}{2 a_n} \sum _{x \in V_n } \varphi(x)  \sum _{y \in V_n}\sum _{z \in V_n:z \sim y}
  \int _0 ^{t}\left( p _n(t-s,x,y) - p_n(t  -s,x, z)\right)
(\h
  _{z} - \h  _y ) (s- ) d A^n _{y,z} (s).
\end{multline}
We remark that due to the graphical representation of the exclusion
process, $f_n$ can be thought of as a function on the probability
space $\left( \{0,1\}^{V_n} \times \O, \cB_n\times \Ff, \nu_n\otimes
\bbP \right)$, where $\Bb_n$ denotes the Borel $\s$--algebra of the
product space  $\{0,1\}^{V_n}$. Moreover, note that  $|f_n|\leq
c(\varphi)$ due to (\ref{formale1}).

In the following  arguments $n$ can be thought of  as fixed. Due to
our assumption on $\nu_n$, given  $\e$ with $0<\e<1$ there exists a
bounded set $C_n \subset X$ such that $\nu _n (A^c)\leq \e$ where
$$
A  = \left\{ \h\,:\, \h_{x} =0\text{ if } x \in V_n \setminus  C_n
\right\}.
$$
Moreover, one can find a bounded set $D_n \subset X$ such that
 $\bbP(\cB ^c) \leq \e$, where
$$
\cB =\left\{\o\,:\, \cup _{x \in V_n \cap  C_n} \cC^{(n)}_x
(t)[\o]\subset D_n \right\}.
$$
 Then $(\nu_n \otimes P) (A\times \cB) \geq (1-\e)^2$. Due to the graphical
representation, one gets  $ \II_{A \times \cB } f_n = \II _{A \times
\cB } \tilde f $,  where
\begin{multline*}
\tilde f= \\
\frac{1}{2 a_n} \sum _{x\in V_n } \varphi(x) \sum _{y \in V_n \cap
C_n}\sum _{z \in V_n : z \sim y }  \int _0 ^{t } \Big( p _n(t -s,x,
y) - p_n(t -s,x, z)\bigr) (\h _{z} - \h _y ) (s- ) d A^n _{y,z} (s).
\end{multline*}
In particular,
$$
\bbE_{\nu_n} (f_n ^2)\leq c(\varphi)^2 \left(\nu _n \times
\bbP\right) ((A \times \cB)^c)+ \bbE_{\nu_n}\left( \II_{A \times
\cB}  \tilde f^2  \right) \leq c(\varphi)^2 (2\e-\e^2)+ \bbE_{\nu_n}
\bigl(\tilde f^2 \bigr).
$$
Since the quadratic variation of the martingale $A^n_b$ is $c(b)t$
and since all series are finite sums, we can compute
\begin{multline*}\bbE_{\nu_n} \bigl(\tilde f^2 \bigr)= \frac{1}{4 a^2_n}\sum _{y
\in V_n \cap C_n}\sum _{z \in V_n : z \sim y } c(y,z) \times \\
 \int
_0 ^{t }ds  (\h _{z} - \h _y ) ^2(s) \Big[ \sum _{x\in V_n }
\varphi(x) \Big(
p _n(t -s,x, y) - p_n(t -s,x, z)\bigr)  \Big]^2\leq \\
\frac{1}{4 a^2_n}\sum _{y \in V_n \cap C_n}\sum _{z \in V_n : z \sim
y } c(y,z)  \int _0 ^{t }ds   \Big[ \sum _{x\in V_n } \varphi(x)
\Big(
p _n(t -s,x, y) - p_n(t -s,x, z)\Big) \Big]^2=\\
\frac{1}{4 a^2_n}\sum _{y \in V_n \cap C_n}\sum _{z \in V_n : z \sim
y } c(y,z)  \int _0 ^{t }ds \bigl( P_{s} ^n \varphi (y)- P_{s} ^n
\varphi (z)\bigr)^2\leq \frac{1}{ 2 a_n } \int _0^ t ds < P_s ^n
\varphi, -\bbL P_s n \varphi>_{m_n}=\\
 -\frac{1}{ 2 a_n } \int _0^ t ds \frac{d}{ds} < P_s ^n
\varphi, P_s n \varphi>_{m_n}=   \frac{1}{2 a_n} \bigl(  < \varphi,
 \varphi>_{m_n}-   < P_t ^n
\varphi, P_t ^n \varphi>_{m_n}\bigr)\leq  \frac{m_n (\varphi^2)}{2
a_n}\,,
\end{multline*}
where $<\cdot, \cdot>_{m_n}$ denotes the scalar product in $L^2 (V_n
, m_n) $ (recall \eqref{vagamente}). Note that in the second
identity we have used the symmetry \eqref{simone} which implies:
$$
\sum _{x\in V_n } \varphi(x) \Big( p _n(t -s,x, y) - p_n(t -s,x,
z)\Big)= P_{t-s} ^n \varphi (y)- P_{t-s} ^n \varphi (z)\,.
$$

The above bound,  the vague convergence  \eqref{vagamente}  and the
assumption $a_n \to \infty$ imply that $\varlimsup
_{n\uparrow\infty} \bbE_{\nu_n} \left( f_n ^2\right) \leq c(\varphi
) ^2 (2\e-\e^2)$. Since $\e$ is arbitrary, we get the thesis.
\end{proof}

\subsection{Conclusion}
We have now all the tools to prove Theorem \ref{annibale}. Indeed,
by Proposition \ref{gaetano}, in order to derive the thesis we only
need to show that, given $\d>0$,
\begin{equation}\label{prof3}
\lim _{n \uparrow \infty} \bbP_{\mu_n} \Big( \Big| \frac{1}{a_n}
\sum _{x\in V_n} \varphi(x)\sum _{y \in V_n } p_n(t,x,y) \h _y (0) -
\int \varphi(u) (P_t \rho_0) (u) m (du) \Big|> \d \Big)=0\,.
\end{equation}
Due to the symmetry \eqref{simone} we can write
$$
 \frac{1}{a_n}
\sum _{x\in V_n} \varphi(x)\sum _{y \in V_n } p_n(t,x,y) \h _y (0)=
\frac{1}{a_n} \sum _{x\in V_n}  \h_x (0) P^n_t \varphi(x)\,.
$$
 By the homogenization
assumption \eqref{hom1}, we only need to prove that
\begin{equation}\label{prof4}
\lim _{n \uparrow \infty} \mu_n  \Big( \Big| \frac{1}{a_n} \sum
_{x\in V_n}  \h_x (0) P_t \varphi(x)  - \int \varphi(u) (P_t \rho_0)
(u) m (du) \Big|> \d/2 \Big)=0\,.
\end{equation}
Since $W$ has symmetric kernel, we conclude that
$$ \int\varphi(u) (P_t \rho_0)
(u) m (du)=  \int P_t \varphi(u) \rho_0  (u) m (du)\,.
$$
Hence, we only need to prove that
\begin{equation}\label{prof5}
\lim _{n \uparrow \infty} \mu_n  \Big( \Big|\frac{1}{a_n} \sum
_{x\in V_n}  \h_x (0) P_t \varphi(x)    - \int \rho_0(u) P_t
\varphi(u)  m (du) \Big|> \d/2 \Big)=0\,.
\end{equation}
This follows from assumption \eqref{prof1}.


\begin{thebibliography}{99}



\bibitem[A]{A}  G.\ Allaire. {\em Homogenization and two--scale
convergence}. SIAM J. Math. Anal., {\bf 23}, 1482--1518 (1992).


\bibitem[D]{D} R.\ Durrett. {\em Ten lectures on particle systems}.
In: {\em Ecole d'Et\'{e} de Probabilit\'{e}s de Saint--Flour
XXIII--1993}, P. Berndard (ed.), Lect. Notes Math. {\bf 1608},
97--201 (1995).



\bibitem[F1]{F1} A.\, Faggionato. \emph{Bulk diffusion of 1D exclusion process with bond
disorder.} Markov Processes and Related Fields {\bf 13}, 519--542
(2007).

\bibitem[F2]{F2} A.\ Faggionato. \emph{Random walks and  exclusion processes among random
conductances on random infinite clusters: homogenization and
hydrodynamic limit.} Electronic Journal of Probability {\bf 13},
2217--2247 (2008).

\bibitem[FJL]{FJL}  A.\ Faggionato, D.M.\ Jara., C.\ Landim \emph{Hydrodynamic behavior
of 1D
 subdiffusive exclusion processes
  with random conductances}.
 Prob. Theory and Related Fields 144, 633-667 (2009).

\bibitem[FM]{FM}
A.\ Faggionato, F.\ Martinelli. \emph{Hydrodynamic limit of a
disordered lattice gas}.  Probab. Theory and Related Fields  {\bf
127} (3), 535--608 (2003).


%\bibitem[Fr]{Fr} J.\  Fritz. {\em  Hydrodynamics in a symmetric random medium}.
%  Comm. Math. Phys.  {\bf 125}  13--25 (1989).

\bibitem[GJ]{GJ}
 P.\ Goncalves, M.\ Jara. \emph{ Scaling limit of gradient systems in random environment}. J. Stat. Phys. {\bf 131}(4),
691--716 (2008).



\bibitem[J]{J} M.\ Jara. {\em Hydrodynamic limit of the exclusion process in
inhomogeneous media.} Preprint, arXiv:0908.4120v1 (2009).

\bibitem[JL]{JL} M. Jara, C. Landim. {\em  Nonequilibrium central limit theorem
  for a tagged particle in symmetric simple exclusion}.
   Annales de l'institut Henri Poincar\'{e} (B) Probabilit\'{e}s et Statistiques, {\bf 42},
   567--577, (2006).



\bibitem[KK]{KK} K.\ Kawazu, H.\ Kesten, {\em On birth and death processes
    in symmetric random environment}. J. Statist. Phys. {\bf 37},
  561--576 (1984).

\bibitem[KL]{KL}
C.\ Kipnis, C.\ Landim. \emph{ Scaling limits of interacting
particle systems}. Springer, Berlin (1999).


\bibitem[L]{L} T.M.\, Liggett. \emph{Interacting particle systems}.
 Springer,  New York  (1985).

\bibitem[M]{M} P.\ Mathieu. \emph{Quenched invariance principles for random walks with random
conductances}.  J. Stat. Phys. {\bf 130}, 1025-–1046 (2008).



\bibitem[N]{N} K.\ Nagy. \emph{Symmetric random walk in random
environment}. Period. Math. Hung. {\bf 45}, 101--120 (2002).

\bibitem[Nu]{Nu} G.\ Nguetseng. \emph{A general convergence result
for a functional related to the theory of homogenization}. SIAM J.
Math. Anal., {\bf 20}, 608--623 (1989).

\bibitem[P]{P} S.E.\ Pastukhova. \emph{On the convergence of
hyperbolic semigroups in a variable Hilbert space}.  J. Math. Sci.
(N.Y.), {\bf 127}, no. 5, 2263--2283 (2005).




\bibitem[Q1]{Q1}
J.\ Quastel. \emph{Diffusion in disordered media}. In
\emph{Proceedings in Nonlinear Stochastic PDEs} (T. Funaki and W.
Woyczinky, eds), Springer, New York, 65--79 (1996).


\bibitem[Q2]{Q2}
J.\ Quastel, \emph{Bulk diffusion in a system with site disorder.}
Ann. Probab. {\bf 34} (5),   1990--2036  (2006)

\bibitem[S]{S} C.\ Stone, {\em Limit theorems for random walks, birth and
    death processes, and diffusion processes}. Ill. J. Math. {\bf 7},
  638--660 (1963).



\bibitem[ZP]{ZP}  V.V.\, Zhikov, A.L.\, Pyatnitskii. \emph{Homogenization
of random singular structures and random measures}. (Russian) Izv.
Ross. Akad. Nauk Ser. Mat. {\bf 70}, no. 1, 23--74 (2006);
translation in Izv. Math. {\bf 70}, no. 1, 19--67 (2006).
\end{thebibliography}
\end{document}